\documentclass[a4paper,12pt]{article}
\usepackage{amssymb}
\def\beq{\begin{equation}}
\def\eeq{\end{equation}}
\def\bea{\begin{eqnarray}}
\def\eea{\end{eqnarray}}
\def\nn{\nonumber}
\def\Uh{U_h(osp(2/1))}
\def\t{{\bf t}}
\def\G{{\cal G}}
\def\JMP{\textit{J. Math. Phys.} }
\def\JPA{\textit{J. Phys. {\bf A}: Math. Gen.} }
\def\LMP{\textit{Lett. Math. Phys.} }
\def\MPL{\textit{Mod. Phys. Lett.} }
\def\IJMP{\textit{Int. J. Mod. Phys.} }
\renewcommand{\theequation}{\thesection.\arabic{equation}}

\begin{document}
\setlength{\baselineskip}{18pt}

\begin{titlepage}
\vspace*{2cm}
\begin{center}
{\Large\bf Super-Jordanian Quantum Superalgebra $ U_h(osp(2/1)) $}

\vspace{1cm}
{\large
N. Aizawa${}^{\dag}$, \ R. Chakrabarti${}^{\ddag}$ and J. Segar${}^{\star}$

\vspace{7mm}
${}^{\dag}$ Department of Applied Mathematics,

Osaka Women's University,
Sakai, Osaka 590-0035, JAPAN

\bigskip
${}^{\ddag}$ Department of Theoretical Physics,
University of Madras,

Guindy Campus,
Chennai, 600 025, INDIA

\bigskip
${}^{\star}$ Department of Physics, R. K. M. Vivekananda College,

Mylapore, Chennai, 600 003, INDIA
}
\end{center}
\vfill
\begin{abstract}
  A triangular quantum deformation of $ osp(2/1) $ from the
classical $r$-matrix including an odd generator is presented with
its full Hopf algebra structure. The deformation maps, twisting
element and tensor operators are considered for the deformed
$ osp(2/1)$. It is also shown that its subalgebra generated by the
Borel subalgebra is self-dual.
\end{abstract}
\end{titlepage}
%
%%%%%%%%%%%%%%%%%%%%%%%%%%%%%%%%%%%%%%%%%%%%%%%%%%%%%%%%%%%%%%%
%
%        Introduction
%
%%%%%%%%%%%%%%%%%%%%%%%%%%%%%%%%%%%%%%%%%%%%%%%%%%%%%%%%%%%%%%%%
%
\setcounter{equation}{0}
\section{Introduction}

  Both Lie and quantum superalgebras play important roles in
various contexts of theoretical physics. These objects are also of
great interest in modern mathematics. The simplest and the most
fundamental Lie superalgebra is $ osp(2/1), $ and it is in a
similar position to $ sl(2) $ in the theory of Lie superalgebra.
Namely, the understanding of its structure and representations is
the basics to the higher rank superalgebras. This means that the
quantization of $ osp(2/1) $ is also fundamental to study quantum
superalgebras. The recent work shows that there exist three
distinct bialgebra structures on $ osp(2/1) $ and all of them are
coboundary \cite{JS}. We therefore have three distinct quantization
for $osp(2/1)$. Two of them are wellknown \cite{KR,CK} but one is
still incomplete. Despite a lot of investigations of quantum
superalgebras, it is a bit surprising that we do not have
a complete list of quantization of $osp(2/1).$  In this paper, we
wish to complete the list by giving a full Hopf superalgebra
structure of the not yet completed quantization of $ osp(2/1). $
We also discuss some algebraic properties of this quantization of
$ osp(2/1)$ such as twisting element, tensor operators,
self-duality and so on.

  The plan of this paper is as follows. The next section is a short review
of bialgebra structures on $ osp(2/1) $ and its quantization. We shall make
clear known and unknown results on the quantization.
In \S 3, a quantization of $ osp(2/1) $ is presented with deformed commutation
relations and a full Hopf structure. This completes the list of quantization of
$ osp(2/1). $
The quantized $ osp(2/1) $ is triangular
so that it has a
basis satisfying undeformed commutation relations.
The maps connecting the bases of deformed and undeformed commutation relations
are given in \S 4. The twisting element for this quantization is investigated
and given as power series in the deformation parameter up to order three in
\S 5. In \S 6, the tensor operators for adjoint representation of the quantized
$ osp(2/1) $ are given. In \S 7, we discuss self-duality of the Borel subalgebra
of the quantized $ osp(2/1). $

%
%%%%%%%%%%%%%%%%%%%%%%%%%%%%%%%%%%%%%%%%%%%%%%%%%%%%%%%%%%%%%%%
%
%        Bialgebras on osp(2/1) and quantization
%
%%%%%%%%%%%%%%%%%%%%%%%%%%%%%%%%%%%%%%%%%%%%%%%%%%%%%%%%%%%%%%%%
%
\setcounter{equation}{0}
\section{Bialgebras on $osp(2/1)$ and quantization}

  The Lie superalgebra $ osp(2/1) $ has three even and two odd elements.
Let $ J_0,\; J_{\pm} $ be even elements and $ v_{\pm} $ be odd ones, then
$ osp(2/1) $ is defined by the relations
\bea
  & & [J_0, v_{\pm}] = \pm \frac{1}{2} v_{\pm}, \qquad
      \{ v_+, v_- \} = -\frac{1}{2} J_0,
  \nn \\
  & & [J_0, J_{\pm} ] = \pm J_{\pm}, \qquad [J_+, J_-] = 2J_0,
  \nn \\
  & & J_{\pm} = \pm 4 v_{\pm}^2,\qquad
      [J_{\pm}, v_{\pm} ] = 0, \qquad
      [J_{\pm}, v_{\mp} ] = v_{\pm}.  \label{udcom}
\eea
The even elements span the $sl(2)$ subalgebra.
Note that $ osp(2/1) $ and its universal enveloping algebra are
generated by $ v_{\pm} $ and $J_0.$

 The bialgebra structures on $osp(2/1)$ are classified in \cite{JS}.
The authors, by using computer, show that all the possible bialgebra
structures are coboundary and there are three inequivalent classical $r$-matrices.
\bea
  & & r_1 = J_0 \wedge J_+, \nn \\
  & & r_2 = J_0 \wedge J_+ - v_+ \wedge v_+, \nn \\
  & & r_3 = t(J_0 \wedge J_+ + J_0 \wedge J_- - v_+ \wedge v_+ - v_- \wedge v_-).
      \label{r123}
\eea
The parameter $t$ in $r_3$ becomes irrelevant in quantization, since it can be
absorbed into a deformation parameter. However, it cannot be removed by a change of
the basis in $ osp(2/1).$  The $r$-matrices $ r_1, r_2 $ satisfy the classical Yang-Baxter
equation, while $ r_3 $ satisfies the modified classical Yang-Baxter equation.
The quantization of these bialgebras has been considered. The $q$-deformation of
$ osp(2/1) $ considered in \cite{KR} corresponds to $r_3$ and it is a quasi-triangular
Hopf superalgebra. The irreducible representations and some applications of the $q$-deformed
$ osp(2/1) $ are studied in \cite{Sal}.
The $r$-matrix $r_1$ consists of the elements in $ sl(2) $ subalgebra and $r_1$
is identical to the classical $r$-matrix of the Jordanian deformation of $ sl(2). $
Thus one can quantize $osp(2/1)$ using the inclusion $ sl(2) \subset osp(2/1). $
This is done in \cite{CK} where the Drinfel'd twist\cite{dri} for $ sl(2) $ is applied to $osp(2/1)$ and
the obtained Hopf superalgebra is triangular.
For these two quantum superalgebras, their full Hopf structures, universal $R$-matrices and
dual quantum supergroups have been obtained.

  On the other hand, the quantization of $r_2$ remains unfinished, though there
are two publications \cite{Kul,JS2}.
Kulish tried to find a twisting element that gives rise to the universal $R$-marix whose
classical counterpart is $ r_2$ \cite{Kul}.
He observed
that the quantum $R$-matrix in fundamental representation corresponding to $r_2$
can be obtained from the one for $q$-deformed $osp(2/1)$ by contraction.
This $R$-matrix is triangular so that consistent with that $r_2$ solves
the classical Yang-Baxter equation. With this $R$-matrix and the FRT-method \cite{FRT},
the subalgebra generated by $ J_0 $ and $v_+$ is quantized.
It is also shown that this $R$-matrix
can be decomposed into Jordanian $sl(2)$ part and additional factor. This suggest
that the twisting element would be a product of $sl(2)$ Jordanian twist and
a factor containing $v_+$. The conjectured form of the twisting element is
given in \cite{Kul}, but it contains unknown factors.
Juszczak and Sobczyk take the same $R$-matrix to construct the dual quantum
supergroup \cite{JS2}.
By the FRT-method, the quantum supergroup was explicitly constructed.
Then from the duality, they obtain the quantization of the subalgebra generated
by $ J_0, J_+ $ and $v_+$. In the next section, we shall give full quantization
of $osp(2/1)$.

%
%%%%%%%%%%%%%%%%%%%%%%%%%%%%%%%%%%%%%%%%%%%%%%%%%%%%%%%%%%%%%%%
%
%        Quantum superalgebra Uh(osp(2/1))
%
%%%%%%%%%%%%%%%%%%%%%%%%%%%%%%%%%%%%%%%%%%%%%%%%%%%%%%%%%%%%%%%%
%
\setcounter{equation}{0}
\section{Quantum superalgebra $U_h(osp(2/1))$}

  We denote the generators of the quantized $ osp(2/1) $ by
$(H, V_{\pm})$. Their classical limits are given by $(J_0, v_{\pm})$,
respectively. We introduce two additional elements $(X, Y)$
whose classical limit are $(J_+ , J_-), $ respectively. With the same
notations, the classical
bialgebra structure is summarized as follows: The classical $r$-matrix reads
\beq
  r = H \otimes X - X \otimes H - 2 V_{+} \otimes V_{+}.  \label{cr}
\eeq
The classical co-commutators, defined as
\beq
  \delta (g) = [r, (g \otimes 1 + 1 \otimes g)], \qquad g \in osp(2/1), \label{cocom}
\eeq
are given by
\bea
 & & \delta (H) = - r, \nn \\
 & & \delta (X) = 0, \nn \\
 & & \delta (Y) =  X \otimes Y - Y \otimes  X + 2 V_{+} \otimes V_{-}
           + 2 V_{-} \otimes V_{+}, \nn \\
 & & \delta (V_{+}) = \frac{1}{2} X \otimes V_{+} - \frac{1}{2} V_{+} \otimes X, \nn \\
 & & \delta (V_{-}) = \frac{1}{2} X \otimes V_{-} - \frac{1}{2} V_{-} \otimes X.
  \label{sdel}
\eea

  We take a quite naive approach for quantization. We have all classical data of
$osp(2/1)$ bialgebra. The quantized algebra must have the
following properties:
(i) The classical co-commutator must be maintained in the classical limit.
(ii) The quantum coproduct is a homomorphism of the algebra, and
coassociative.
(iii) We know from the previous works \cite{Kul,JS2} that
 $X$ must be primitive.
With these informations, one can write down the algebraic relations and the
quantum coproducts. A similar approach for quantization is discussed in \cite{LM}
where quantum coproducts obtained first then commutation relations are derived.
In the present case, the commutation relations read
\bea
 & & [H, V_{+}] = \frac{1}{2} V_{+} \cosh (hX),\qquad
     [H, V_{-}] = - \frac{1}{4} V_{-} \cosh (hX) - \frac{1}{4} \cosh (hX) V_{-},\nn \\
 & & \{ V_{+}, V_{-} \} = - \frac{1}{2} H, \qquad\qquad\quad
     [H, X] = \frac{1}{h} \sinh(hX), \nn \\
 & & [H, Y] = - \frac{1}{2} Y \cosh (hX) - \frac{1}{2} \cosh (hX) Y
         + h V_{-} \sinh (hX) V_{+}
         - h V_{+} \sinh (hX) V_{-},
  \nn \\
 & & [X, Y] = 2 H,   \qquad\qquad
     V_{+}^{2} = \frac{1}{4h} \sinh (hX), \qquad\qquad
     V_{-}^{2} = - \frac{1}{4} Y, \nn \\
 & & [X, V_{+}] = [Y, V_{-}] = 0,  \quad
     [X, V_{-}] = V_{+}, \nn \\
 & & [Y, V_{+}] = \frac{1}{2} V_{-} \cosh (hX) + \frac{1}{2} \cosh (hX) V_{-},
     \label{alg}
\eea
where $h$ is the deformation parameter and $ h \rightarrow 0 $ gives the
classical limit. The quantum coproducts $\Delta$ read
\bea
 & & \Delta (H) = H \otimes T^{-1} + T \otimes H
             + 2h V_{+} T^{1/2} \otimes V_{+} T^{-1/2},
 \nn \\
 & & \Delta (X) = X \otimes 1 + 1 \otimes X,
 \nn \\
 & & \Delta (Y) = Y \otimes T^{-1} + T \otimes Y
             + 2h V_{+} T^{1/2} \otimes T^{-1/2} V_{-}
   + 2h T^{1/2} V_{-} \otimes V_{+} T^{-1/2},
  \nn \\
 & & \Delta (V_{\pm}) = V_{\pm} \otimes T^{-1/2} + T^{1/2} \otimes V_{\pm},
  \label{copro}
\eea
where $ T = \exp(hX).$
We see that there is only one primitive element.
It is straightforward to verify that the commutation relations and the
coproduct are consistent with the axioms.
The counit $\epsilon$ and the antipode $S$ follow form the coproduct
and they are given by
\beq
  \epsilon(H) = \epsilon(V_{\pm}) = \epsilon(X) = \epsilon(Y) = 0,
  \label{counit}
\eeq
and
\beq
  \begin{array}{lclcl}
      S(H) = -H - 2hV_+^2, & \ \ &
      S(V_+) = -V_+,       & \ &
      S(V_-) = -V_- + \frac{h}{2}V_+ \\
      S(X) = -X,           & &
      S(Y) = -Y + hH + h^2 V_+^2. & &
  \end{array}
  \label{antipode}
\eeq
The relations (\ref{copro})-(\ref{antipode}) define a quantum superalgebra.
We denote this algebra by $ \Uh. $ By construction, the algebra $\Uh$
is a triangular Hopf superalgebra and a one parameter deformation of
the algebra $ U(ops(2/1))$ so that the algebra $\Uh $ has three generators,
$ H $ and $ V_{\pm}. $

%
%%%%%%%%%%%%%%%%%%%%%%%%%%%%%%%%%%%%%%%%%%%%%%%%%%%%%%%%%%%%%%%
%
%        Deformation maps
%
%%%%%%%%%%%%%%%%%%%%%%%%%%%%%%%%%%%%%%%%%%%%%%%%%%%%%%%%%%%%%%%%
%
\setcounter{equation}{0}
\section{Deformation maps}

  It is known that a triangular Hopf algebra can be obtained from
a Lie algebra by Drinfel'd twist \cite{dri}. The algebra obtained by
Drinfel'd twist has deformed coproduct and deformed antipode, while
the commutation relations and the counit remain undeformed.
This can be extended to superalgebras. The quantum superalgebra $ \Uh $
defined above is a triangular Hopf superalgebra and has deformed
commutation relations. This means that there exist maps that connect
generators satisfying deformed commutation relations and undeformed
ones. In this section, we seek such deformation maps. We first give
a general class of such maps, then give two explicit examples.

An ansatz for a general class of maps may be assumed as
\begin{eqnarray}
&&V_{+} = f_{1} (J_{+})\, v_{+}\qquad H = f_{2} (J_{+})\, J_{0},
\nonumber\\
&&V_{-} = f_{3} (J_{+})\, v_{-} + u(J_{+})\, v_{+} + w (J_{+})\,
v_{+} J_{0},
\label{eq:genmap}
\end{eqnarray}
where we introduce $(f_{1}, f_{2}, f_{3}; u, w)$ as functions of
$J_{+}$ only. The elements of $ \Uh $ with undeformed commutation
relations are denoted by the same notations as the classical
$ osp(2/1)$. An additive function $f(J_{+})$ in the expression of
$H$ may be absorbed by a similarity transformation. To ensure correct
classical limits, the above introduced functions are  required to
satisfy the limiting properties:
\begin{equation}
(f_{1}, f_{2}, f_{3}; u, w)\,\rightarrow (1, 1, 1; 0, 0)
\label{eq:flimit}
\end{equation}
as $h \rightarrow 0$. The operators $T^{\pm 1} \equiv \exp(\pm hX)$
may now be expressed as
\begin{eqnarray}
T &=& h J_{+}\,(f_{1}(J_{+}))^{2} + \sqrt{1 + h^2 J_{+}^{2}\,
(f_{1}(J_{+}))^{4}},\nonumber\\
T^{-1} &=& - h J_{+}\,(f_{1}(J_{+}))^{2} + \sqrt{1 + h^2 J_{+}^{2}\,
(f_{1}(J_{+}))^{4}}.
\label{eq:TTinv}
\end{eqnarray}
Substituting the ansatz (\ref{eq:genmap}) in the defining algebraic
relations systematically, we, for a {\it given} function $f_{1}$,
obtain a set of {\it six} nonlinear equations for {\it four} unknown
functions $(f_{2}, f_{3}; u, w)$:
\begin{eqnarray}
&&f_{2}(J_{+}) \left(2 J_{+}\,f_{1}^{\prime}(J_{+}) +
f_{1}(J_{+})\right) - \sqrt{ 1 + h^{2}J_{+}^{2} (f_{1}(J_{+}))^{4}}\,
f_{1}(J_{+}) = 0,\nonumber\\
&&2 J_{+}\,f_{2}(J_{+})\,f_{3}^{\prime} (J_{+}) - f_{2}(J_{+})\,
f_{3}(J_{+}) + \sqrt{1 + h^{2} J_{+}^{2} (f_{1}(J_{+}))^{4}}\,
f_{3}(J_{+}) = 0,\nonumber\\
&&2 J_{+}\,f_{2}(J_{+})\,u^{\prime}(J_{+}) + f_{2}(J_{+})\,u(J_{+})
+ \sqrt{1 + h^{2} J_{+}^{2} (f_{1}(J_{+}))^{4}}\,u(J_{+})\nonumber\\
&&\phantom{qqqqqqqqqqqqqqqqqqq}-
\frac{h^{2}}{2}J_{+}\,(f_{1}(J_{+}))^{3} = 0,\nonumber\\
&&J_{+}\,f_{1}(J_{+})\,w(J_{+}) - f_{1}(J_{+})\,f_{3}(J_{+})
+ f_{2}(J_{+}) = 0,\nonumber\\
&&2 J_{+}\,f_{1}(J_{+})\,u(J_{+}) + J_{+}^{2}\,
f_{1}^{\prime}(J_{+})\,w(J_{+}) - J_{+}\,f_{1}^{\prime}(J_{+})\,
f_{3}(J_{+})\nonumber\\
&&\phantom{qqqqqqqqqqqqqqqqqqq}+ \frac{1}{2} J_{+}\,
f_{1}(J_{+})\,w(J_{+}) =0,\nonumber\\
&&f_{3}(J_{+})\,f_{2}^{\prime}(J_{+}) + J_{+}\,f_{2}(J_{+})\,
w^{\prime}(J_{+}) - J_{+}\,w(J_{+})\,f_{2}^{\prime}(J_{+})
+ \frac{1}{2} f_{2}(J_{+})\,w(J_{+}) \nonumber\\
&&\phantom{qqqqqqqqqqqqqqqqqqq}+ \frac{1}{2}\,\sqrt{1 + h^{2}
J_{+}^{2}\,(f_{1}(J_{+}))^{4}}\,w(J_{+}) = 0.
\label{eq:fdiff}
\end{eqnarray}
We made an identification $ {\displaystyle [J_0, f(J_+)] = J_+ \frac{d}{dJ_+}f(J_+)} $
on a function $f(J_+)$.
These equations may then be solved consistently and the classical
limit (\ref{eq:flimit}) may be implemented through proper choice
of integration constants. This leads to the following solution:
\begin{eqnarray}
&&f_{2} (J_{+}) = \frac{\sqrt{1 + h^{2} J_{+}^{2}\,
(f_{1}(J_{+}))^{4}}}
{f_{1}(J_{+}) + 2 J_{+}\,f_{1}^{\prime}(J_{+})}\,f_{1}(J_{+}),
\qquad f_{3}(J_{+}) = \frac{1}{f_{1}(J_{+})},\nonumber\\
&&u(J_{+}) = -\frac{1}{4}\,w(J_{+}) + \frac{1}{2}\,
\frac{f_{1}^{\prime}(J_{+})}{(f_{1}(J_{+}))^{2}}\,f_{2}(J_{+}),
\qquad w(J_{+}) = \frac{1 - f_{2}(J_{+})}{J_{+}\,f_{1}(J_{+})}.
\label{eq:fsolve}
\end{eqnarray}
We see that for a given function $ f_1$ other functions are uniquely determined.
This shows that deformation maps are not unique.

\par

Even for the invertible maps it is useful to consider the general
solution starting from the other end. So we consider the following
class of inverse maps given by
\begin{eqnarray}
&&v_{+} = g_{1} (T)\, V_{+},\qquad J_{0} = g_{2} (T)\, H,\nonumber\\
&&v_{-} = g_{3} (T)\, V_{-} + a (T)\, V_{+} + b (T)\, V_{+} H,
\label{eq:ingenmap}
\end{eqnarray}
where $(g_{1}, g_{2}, g_{3}; a, b)$ are functions of $T$ only. In the
classical limit $h \rightarrow 0$, these functions obey the
following limiting properties
\begin{equation}
(g_{1}, g_{2}, g_{3}; a, b)\,\rightarrow (1, 1, 1; 0, 0).
\label{eq:glimit}
\end{equation}
Substituting the above ansatz (\ref{eq:ingenmap}) in the defining
algebra, we, as before, obtain a set of six coupled nonlinear
equations:
\begin{eqnarray}
&&(T^{2} - 1)\,g_{1}^{\prime}(T)\,g_{2}(T) + \frac{1}{2}\,
(T + T^{-1})\,g_{1}(T)\,g_{2}(T) - g_{1}(T) =0,\nonumber\\
&&(T^{2} - 1)\,g_{2}(T)\,g_{3}^{\prime}(T) - \frac{1}{2}\,
(T + T^{-1})\,g_{2}(T)\,g_{3}(T) + g_{3}(T) =0,\nonumber\\
&&(T^{2} - 1)\,g_{2}(T)\,a^{\prime}(T) + \left(1 + \frac{1}{2}\,
(T + T^{-1})\,g_{2}(T)\right)\, a(T) + \frac{h}{4}
(T - T^{-1})\,g_{2}(T)\,g_{3}(T) =0,\nonumber\\
&&(T^{2} - 1)\,g_{2}(T)\,b^{\prime}(T) + \left(1 + \frac{1}{2}\,
(T + T^{-1})\,g_{2} (T)\,- (T^{2} - 1)\, g_{2}^{\prime}(T) \right)\,
b(T) + 2 hT\,g_{2}^{\prime}(T)\,g_{3}(T) = 0,\nonumber\\
&&(T - T^{-1})\,g_{1}(T)\,b(T) + 2 h ( g_{2}(T) - 1 ) = 0,
\nonumber\\
&&\left((T + T^{-1})\,g_{1}(T) + 2 (T^{2} - 1)\,g_{1}^{\prime}(T)
\right)\,b(T) + 8 g_{1}(T)\,a(T) - 4 hT\,g_{1}^{\prime}(T)\,
g_{3}(T) = 0,
\label{eq:gset}
\end{eqnarray}
where we have used the identity
$ {\displaystyle [H, g(T)] = \frac{1}{2}(T^2-1) \frac{d}{dT}g(T)}$.
Treating the function $g_{1}(T)$ as known and taking the boundary
conditions (\ref{eq:glimit}) into account the remaining functions
may be solved unambiguously as follows:
\begin{eqnarray}
&&g_{2}(T) = 2\,\frac{g_{1}(T)}{(T+T^{-1})\,g_{1}(T) + 2 (T^{2} - 1)
\,g_{1}^{\prime}(T)}\qquad g_{3}(T) = \frac{1}{g_{1}(T)}\nonumber\\
&&a(T) = - \frac{h}{4}\,\tanh\left(\frac{h X}{2}\right)\,
(g_{1}(T))^{-1},\qquad b(T) = \frac{2 h}{T - T^{-1}}\,
\frac{1 - g_{2}(T)}{g_{1}(T)}.
\label{eq:gsolve}
\end{eqnarray}

We now give two examples of the maps obtained above.
The first choice of $f_1(J_+)$ and $ g_1(T) $ is
\beq
  f_{1}(J_{+}) = \frac{1}{(1 - 2h J_{+})^{1/4}},\qquad
  g_{1}(T) = T^{-1/2},
  \label{1stc}
\eeq
and the second choice is
\beq
  f_{1}(J_{+}) = \frac{1}{\sqrt{1 - \frac{h^{2}}{4} J_{+}^{2}}}, \qquad
  g_{1}(T) = \hbox{sech}\left(\frac{h X}{2}\right).
  \label{2ndc}
\eeq
The first choice gives the map
\bea
  & & V_+ = e^{\sigma/4} v_+, \qquad H = e^{-\sigma/2} J_0,
  \nn \\
  & & V_- = e^{-\sigma/4} v_- + \frac{h}{4}e^{\sigma/4}
          \tanh\left(\frac{\sigma}{4}\right) v_+
          + \frac{h}{\cosh(\sigma/4)} v_+ J_0,
  \nn \\
  & & X = \frac{1}{2h} \sigma, \qquad
      Y = -4V_-^2, \qquad \sigma = - \ln(1-2hJ_+),
  \label{defmap1}
\eea
and its inverse
\bea
  & & v_+ = e^{-hX/2} V_+, \qquad J_0 = e^{hX}{H},\nn \\
  & & v_- = e^{hX/2} V_- - \frac{h}{4} e^{hX/2} \tanh\left(\frac{hX}{2}\right)
            V_+ - \frac{2he^{hX/2}}{1+e^{-hX}} V_+ H,
      \label{inmap} \\
  & & J_{\pm} = \pm 4 v_{\pm}^2. \nn
\eea
We call the map (\ref{defmap1}) first deformation map.
The induced quantum coproducts for  $ J_0, J_{\pm} $
and $ v_{\pm} $ are obtained from (\ref{inmap}) and (\ref{copro}).
\bea
  & & \Delta(J_0) = J_0 \otimes 1 + e^{\sigma} \otimes J_0 +
      2h v_+ e^{\sigma} \otimes v_+ e^{\sigma/2}, \nn \\
  & & \Delta(J_+) = J_+ \otimes e^{-\sigma} + 1 \otimes J_+, \nn \\
  & & \Delta(v_+) = v_+ \otimes e^{-\sigma/2} + 1 \otimes v_+, \nn \\
  & & \Delta(\sigma) = \sigma \otimes 1 + 1 \otimes \sigma. \label{copro2}
\eea
$ \Delta(v_-), \Delta(J_-) $ are quite messy, so we do not give them here.
Now $ \sigma $ is the only primitive element. One can see that, with a slight
change of conventions,  the coproducts (\ref{copro2}) are identical to the
ones in \cite{Kul} and \cite{JS2}. Thus the first deformation map connects the
basis of $\Uh$ in \cite{Kul,JS2} to ours. The map is also a natural extension of
the similar map \cite{NA1} for Jordanian $ sl(2) $ to $ \Uh. $

  The second choice (\ref{2ndc}) gives the map
\bea
  & & V_+ = \frac{1}{\sqrt{1 - \frac{h^{2}}{4} J_{+}^{2}}} v_+, \quad
      H = J_0, \qquad
      V_- = \sqrt{1- ({\textstyle \frac{h}{2}J_+})^2} v_- +
          \frac{h}{4}\frac{\textstyle \frac{h}{2} J_+}{
            \sqrt{1- ({\textstyle \frac{h}{2} J_+})^2} }v_+,
      \nn \\
  & & X = \frac{2}{h} {\rm arctanh}\left(\frac{h}{2}J_+\right)
        = \frac{1}{h} \ln\left(
        \frac{1 + {\textstyle \frac{h}{2} J_+}}{1 - {\textstyle \frac{h}{2} J_+}}
        \right)_,
      \qquad Y = -4V_-^2, \label{anodefm}
\eea
and its inverse
\bea
 & & v_+ = \hbox{sech}\left(\frac{h X}{2}\right) V_+, \qquad
     J_0 = H, \qquad
     v_- = \left(\cosh\frac{hX}{2}\right) V_- -
         \frac{h}{4}\left( \sinh\frac{hX}{2} \right) V_+,
     \nn \\
 & & J_+ = \frac{2}{h} \tanh\frac{hX}{2}, \qquad
     J_- = -4v_-^2. \label{inanodefm}
\eea
We call the map (\ref{anodefm}) second deformation map.
The map is a natural extension of
the similar map \cite{ACC} for Jordanian $ sl(2) $ to $ \Uh. $

  The irreducible representations of $ \Uh $ are identical to the ones
of the classical $ osp(2/1),$
since $ \Uh $ has undeformed commutation relations.
The representation matrices for $ H, X, Y $ and $V_{\pm}$ are
obtained by the deformation maps. Because deformation maps are not
unique, we have many different matrices for a element of $ \Uh. $
As examples,
the fundamental and the adjoint representation matrices by the first
and the second deformation maps are given in Appendix.

%
%%%%%%%%%%%%%%%%%%%%%%%%%%%%%%%%%%%%%%%%%%%%%%%%%%%%%%%%%%%%%%%
%
%        Twisting elements for the maps
%
%%%%%%%%%%%%%%%%%%%%%%%%%%%%%%%%%%%%%%%%%%%%%%%%%%%%%%%%%%%%%%%%
%
\setcounter{equation}{0}
\section{Twisting elements for the maps}

  Twisting elements are the most fundamental objects for triangular
quantum (super) algebras, since all such deformations are caused by
appropriate twists. In other words,
the twisting element contains all information about the deformation.
The universal $R$-matrix is immediately constructed from the twisting
element. In this section, the twisting element for $ \Uh $ is investigated
and it is given as a power series in the deformation parameter $h$.
The series up to $O(h^3)$ is explicitly given. The higher order terms
can, in principle, be obtained in a similar way.
In many literatures, twisting elements are given in terms of the
classical generators. Here the twisting element for $\Uh$ is
presented in terms of deformed generators ($H, X, Y, V_{\pm}$), since,
as seen in the previous section,
the coproducts for deformed generators are far simpler in form than the
undeformed ones. Thus the twisting element has deformation map dependence in its
expression.
We shall give the twisting element for the first and second deformation maps,
respectively. The relation of deformation maps to twisting elements is studied
in \cite{ACCS,CQ} for Jordanian deformation of $ sl(2) $ and $ gl(2).$

  The general structure of the twisting element corresponding to the
given maps may be described as follows. Let $m$ be a deformation map and
$m^{-1}$ be its inverse
  \begin{equation}
  m:(V_{\pm}, H) \rightarrow (v_{\pm}, J_{0}),\qquad
  m^{-1}:(v_{\pm}, J_{0}) \rightarrow (V_{\pm}, H)
  \label{eq:mminv}
\end{equation}
the classical $(\Delta_0)$ cocommutative and the quantum $(\Delta)$
non-cocommutative coproducts may be related {\it via} the twisting
element $ {\cal G}$ as
\begin{equation}
 \G\,\Delta \circ m^{-1}(\phi)\,\G^{-1} = (m^{-1} \otimes m^{-1})
 \circ \Delta_0(\phi)\qquad \forall \phi \in U(osp(2/1),
 \label{eq:maptwist}
\end{equation}
where the twisting element $\G$ satisfies the cocycle condition
\begin{equation}
  (\G \otimes 1)\,((\Delta \otimes \hbox{id})\G)
  = (1 \otimes \G)\,((\hbox{id} \otimes \Delta)\G).
  \label{eq:cocycle}
\end{equation}
In the present case we, corresponding to the maps discussed in the previous
section,
obtain a series expansion of the twisting element in powers of the
deformation parameter $h$:
\begin{equation}
 \G = 1 \otimes 1 + h \G_{1} + h^{2} \G_{2} + h^{3} \G_{3} +\cdots.
\label{eq:series}
\end{equation}
For the first map these expansion coefficients of the twisting element
reads
\begin{eqnarray}
  && \G_{1} = 2 X \otimes H + 2 V_{+} \otimes V_{+},\nonumber\\
  && \G_{2} = \frac{\G_{1}^{2}}{2!} + 2 X \otimes XH + XV_{+}
  \otimes V_{+},\nonumber\\
  && \G_{3} = \frac{\G_{1}^{3}}{3!} +
  \frac{1}{2} \G_{1} (2X \otimes XH + XV_{+} \otimes V_{+}) +
  \frac{1}{2} (2X \otimes XH + XV_{+} \otimes V_{+}) \G_{1}\nonumber\\
  &&\phantom{G_{3} =} + X \otimes X^{2}H - \frac{1}{4} V_{+} \otimes
  X^{2}V_{+} - \frac{1}{12} X^{2}V_{+} \otimes V_{+}
  + \frac{5}{6} XV_{+} \otimes XV_{+}
  \label{eq:twist1}
\end{eqnarray}
and for the second map they are given by
\begin{eqnarray}
  && \G_{1} = X \otimes H - H \otimes X + 2 V_{+} \otimes V_{+},
  \nonumber\\
  && \G_{2} = \frac{\G_{1}^{2}}{2!} + \frac{1}{4} (H \otimes X^{2} +
  X^{2} \otimes H + 2 XV_{+} \otimes V_{+} - 2 V_{+} \otimes XV_{+}),
  \nonumber\\
  && \G_{3} = \frac{\G_{1}^{3}}{3!} +
  \frac{1}{8} \G_{1} (H \otimes X^{2} + X^{2} \otimes H
  + 2 XV_{+} \otimes V_{+} - 2 V_{+} \otimes XV_{+})\nonumber\\
  &&\phantom{\G_{3} =} + \frac{1}{8} (H \otimes X^{2} + X^{2} \otimes H
  + 2 XV_{+} \otimes V_{+} - 2 V_{+} \otimes XV_{+}) \G_{1}\nonumber\\
  &&\phantom{G_{3} =} - \frac{1}{24} (2 XH \otimes X^{2}
  - 2 X^{2} \otimes XH + X^{2}H \otimes X - X \otimes X^{2}H
  - 6 XV_{+} \otimes XV_{+})\nonumber\\
  &&\phantom{\G_{3} =} - \frac{1}{12} (V_{+} \otimes X^{2}V_{+}
  + X^{2}V_{+} \otimes V_{+} + 2 XV_{+} \otimes XV_{+}).
  \label{eq:twist2}
\end{eqnarray}

It may be checked that the twist operators corresponding to the two
maps described above satisfy the cocycle condition (\ref{eq:cocycle})
upto the desired order in $h$. Moreover, following \cite{ACCS},
we may interrelate the above two twist operators pertaining to
two different maps.

\par

Let us, for the purpose of avoiding confusion, denote two mapping
functions, given in (\ref{1stc}) and (\ref{2ndc}), by $g_{1}$
and ${\hat g}_{1}$; and the corresponding twist operators, given in
(\ref{eq:twist1}) and (\ref{eq:twist2}), by $\G$ and $\hat{\G}$,
respectively. Following \cite{ACCS} these mapping functions may be
related by a similarity relation:
\begin{equation}
U^{-1} g_{1}(T) U = \hat{g}_{1}(T)\quad\Rightarrow\quad
U^{-1} T U = f_{1}(\hat{g}_{1}(T)).
\label{eq:simrel}
\end{equation}
As discussed in \cite{ACCS}, for the purpose of demonostrating the
equivalence of two maps, it is {\it sufficient} to ensure the
relation (\ref{eq:simrel}). The transforming operator $U$ may
be expressed in a series by direct computation:
\beq
U = \exp \left( \left( - h X - \frac{1}{4} (h X)^{2}
+ \frac{5}{24} (h X)^{3} \right) H \right).
\label{eq:Uval}
\eeq
Now the two relevant twist operators, described above, may be related
{\it{\` a} la} \cite{ACCS} as
\beq
\hat{\G}  = ( U^{-1} \otimes U^{-1})\, \G \, ( \Delta( U ) ).
\label{eq:GGrel}
\eeq
%%%%%%%%%%%%%%%%%%%%%%%%%%%%%%%%%%%%%%%%%%%%%%%%%%%%%%%%%%%%%%%%%%%%
%
%        Tensor operators
%
%%%%%%%%%%%%%%%%%%%%%%%%%%%%%%%%%%%%%%%%%%%%%%%%%%%%%%%%%%%%%%%%
%
\setcounter{equation}{0}
\section{Tensor operators}

 We obtained the twisting element $ {\cal G}$ up to $O(h^3)$ in the
previous section. The higher order computation
seems to be difficult. As mentioned in the previous section, it is
important to obtain a closed form of the twisting element.
Let us now recall that tensor operators for twisted algebras are obtained from
the undeformed ones by using the twisting element \cite{Fio}.
Let ${\bf t}, \ {\bf t}_0$ be
tensor operators for deformed and undeformed algebras, respectively.
Then
\beq
  {\bf t} = \mu (id \otimes S)({\cal F} ({\bf t}_0 \otimes 1){\cal F}^{-1}),
  \label{tt0}
\eeq
where ${\cal F} = {\cal G}^{-1} $ is the twisting element and
$ \mu $ is the product of algebra,	 $ \mu(a\otimes b) = ab. $
Thus obtaining the explicit form of tensor operators could be the first step
to obtain a closed form of twisting elements.
Note also that tensor operators for various Lie algebras
appear in many context of physics. The consideration of tensor operators
has physical importance.

  We start with the definition of tensor operators.
Let $a, b$ be elements of $\Uh$ and $ {\bf t},
{\bf s}$ be a operator acting on a given space.
We assume that
$ (a \otimes b) ({\bf t} \otimes {\bf s}) = (-)^{p(b)p({\bf t})} a{\bf t} \otimes b {\bf s}, $
where $p(b)$ and $ p({\bf t}) $ take $ \pm 1 $ depending on even or odd.
The adjoint action of $ a$ on ${\bf t}$ is
defined by
\beq
   {\rm ad}\, a ({\bf t}) = \mu (id \otimes S) \Delta(a) ({\bf t} \otimes 1).
   \label{DEFad}
\eeq
It follows from this definition that
\beq
  {\rm ad}\, ab({\bf t}) = {\rm ad}\, a \circ {\rm ad}\, b({\bf t}). \label{adab}
\eeq
Thus ${\rm ad}$ gives a representation of the algebra
\beq
  {\rm ad}\,[a, b]({\bf t}) = [{\rm ad}\,a, {\rm ad}\, b]({\bf t}),
  \qquad
  {\rm ad}\,\{a, b \}({\bf t}) = \{{\rm ad}\,a, {\rm ad}\, b\}({\bf t}). \label{adrep}
\eeq
Let ${\cal I}$ be a set of indices,
$ \{\; {\bf t}_i,\ i \in {\cal I} \; \}$ be a set of
operators and $ D(a) $ be a representation matrix of $a$. If the operators
$ {\bf t}_i $ form a representation basis of the algebra, namely
\beq
  {\rm ad}\, a({\bf t}_i) = \sum_{j \in {\cal I}} D(a)_{ji} {\bf t}_j, \label{deften}
\eeq
we call $ {\bf t}_i $ tensor operators for the representation $D$.

  The explicit expressions of $ \Uh $ adjoint action for an
even operator ${\bf t}$ are given by
\bea
 & & {\rm ad } X({\bf t}) = [X, {\bf t}], \nn \\
 & & {\rm ad } V_+({\bf t}) = [V_+ T^{-1/2},\; T^{1/2}{\bf t}]T^{1/2}, \nn \\
 & & {\rm ad } H({\bf t}) = [HT^{-1},\; T{\bf t}]T - 2h(T{\bf t}V_+^2 + V_+T^{1/2} {\bf t} V_+T^{1/2}),
     \label{adeven} \\
 & & {\rm ad } V_-({\bf t}) = [V_- T^{-1/2},\; T^{1/2}{\bf t}]T^{1/2} + \frac{h}{2}T^{1/2} {\bf t} V_+,
     \nn \\
 & & {\rm ad } Y({\bf t}) = [YT^{-1},\; T{\bf t}]T + \frac{h}{2}\{HT^{-1},\; T{\bf t}\}T
      -2hT^{1/2}(V_+{\bf t}V_- + V_-{\bf t}V_+)T^{1/2} \nn \\
 & & \qquad\quad\ - \frac{h}{2}{\rm ad}\, H({\bf t}),\nn
\eea
and for a odd operator ${\bf t}$
\bea
 & & {\rm ad } X({\bf t}) = [X, {\bf t}], \nn \\
 & & {\rm ad } V_+({\bf t}) = \{V_+ T^{-1/2},\; T^{1/2}{\bf t}\} T^{1/2}, \nn \\
 & & {\rm ad } H({\bf t}) = [HT^{-1},\; T{\bf t}]T - 2h(T{\bf t}V_+^2 - V_+T^{1/2} {\bf t} V_+T^{1/2}),
     \label{adodd} \\
 & & {\rm ad } V_-({\bf t}) = \{V_- T^{-1/2},\; T^{1/2}{\bf t}\} T^{1/2} - \frac{h}{2}T^{1/2} {\bf t} V_+,
     \nn \\
 & & {\rm ad } Y({\bf t}) = [YT^{-1},\; T{\bf t}]T + \frac{h}{2}\{HT^{-1},\; T{\bf t}\}T
      + 2hT^{1/2}(V_+{\bf t}V_- + V_-{\bf t}V_+)T^{1/2} \nn \\
 & & \qquad\quad\ - \frac{h}{2}{\rm ad}\, H({\bf t}).\nn
\eea
Let us consider the tensor operators for the adjoint representation,
since for the representation
tensor operators are given in terms of the elements of algebra
itself. For comparison, we start with tensor operators for
the classical $ osp(2/1). $
With the index set $ {\cal I} = \{1, 1/2, 0, -1/2, -1 \}, $
it is obvious that the tensor operators for the adjoint representation
(\ref{adjrep}) of $ osp(2/1) $ are given by
\beq
  \t_1 = J_+, \qquad \t_{1/2} = v_+, \qquad \t_0 = J_0, \qquad
  \t_{-1/2} = v_-, \qquad \t_{-1} = J_-. \label{tenundef}
\eeq
Since the first and second deformation maps give the different form of matrices
for $ H, X, Y $ and $ V_{\pm},$ we have two different expressions of tensor
operators for $ \Uh. $ Of course, these two expressions are related each other.
The adjoint representation matrices by the first deformation map are found in
(\ref{1stdefmat}). The explicit forms of the tensor operators become
\bea
  & & \t_1 = \frac{1}{h}T^{-1} \sinh(hX), \nn \\
  & & \t_{1/2} = V_+ T^{-3/2}, \nn \\
  & & \t_0 = HT^{-1} + \frac{3}{4}T^{-1}\sinh(hX), \nn \\
  & & \t_{-1/2} = V_- T^{1/2} + 2hHV_+ T^{-1/2} + \frac{h}{2}V_+ T^{-1/2} \sinh(hX), \nn \\
  & & \t_{-1} = YT + 2hH^2 -2hV_- V_+ T + hH(\cosh(hX) + T^{-1}) \nn \\
  & & \quad \ \  + \frac{h}{8}(5T^{-1} - 3\sinh(hX)) \sinh(hX). \label{tenfor1}
\eea
The tensor operators for the adjointe representation (\ref{2nddefmat}) by the
second deformation map are given by
\bea
  & & \mbox{\boldmath{$\tau_1$}} = \t_1 = \frac{1}{h}T^{-1} \sinh(hX), \nn \\
  & & \mbox{\boldmath{$\tau_{1/2}$}} = \t_{1/2} = V_+ T^{-3/2}, \nn \\
  & & \mbox{\boldmath{$\tau_0$}} = \t_0 = HT^{-1} + \frac{3}{4}T^{-1}\sinh(hX), \nn \\
  & & \mbox{\boldmath{$\tau_{-1/2}$}} = \t_{-1/2} - \frac{h}{2} \t_{1/2} =
       V_- T^{1/2} + 2hHV_+ T^{-1/2} + \frac{h}{2}V_+T^{-1/2}(\sinh(hX) - T^{-1}), \nn \\
  & & \mbox{\boldmath{$\tau_{-1}$}} = \t_{-1} - 2h\t_0 + \frac{h^2}{2} \t_1 \nn \\
  & & \quad \ \ = YT + 2hH^2 -2hV_- V_+ T + hH\sinh(hX) - \frac{3h}{8}\sinh(hX) \cosh(hX).
      \label{tenfor2}
\eea
The relation to (\ref{tenfor1}) is also given in (\ref{tenfor2}).
It can be seen that the first three tensor operators, that reduced to the Borel
subalgebra of $ osp(2/1) $ in the classical limit, are identical in both
expressions.

%
%%%%%%%%%%%%%%%%%%%%%%%%%%%%%%%%%%%%%%%%%%%%%%%%%%%%%%%%%%%%%%%
%
%        Self-duality of Borel subalgebra
%
%%%%%%%%%%%%%%%%%%%%%%%%%%%%%%%%%%%%%%%%%%%%%%%%%%%%%%%%%%%%%%%%
%
\setcounter{equation}{0}
\section{Self-duality of Borel subalgebra}

 It is known that the Jordanian quantization of $ sl(2)$ Borel subalgebra generated by
$ J_0, J_+ $ is self-dual \cite{Vla1}. We show the similar is true for $ \Uh. $
Let $ {\cal B} $ be a quantization of the Borel subalgebra of $ osp(2/1)$
generated by $ H, X $ and $ V_+$
and
$\{\; E_{k\ell m} = H^k X^{\ell} V_+^m,\ k, \ell \in {\mathbb Z}_{\geq 0},\ m = 0, 1 \; \} $
 be its basis.
Let
$\{\; e^{k\ell m} = H^k X^{\ell} V_+^m,\ k, \ell \in {\mathbb Z}_{\geq 0},\ m = 0, 1 \; \} $
be the basis of the algebra $ {\cal B}^* $ dual to $ {\cal B} $ such that
$ <E_{k\ell m}, e^{pqr}> = \delta_{k,p} \delta_{\ell,q} \delta_{m,r}. $
Then by definition of the duality
\bea
 & & E_{k\ell m} E_{k'\ell' m'} = f_{k\ell m,\; k'\ell' m'}^{\quad pqr} E_{pqr},
     \qquad
     \Delta(E_{k\ell m}) = g^{pqr,\; p'q'r'}_{\quad k\ell m} E_{pqr} \otimes E_{p'q'r'},
     \nn \\
 & & e^{pqr} e^{p'q'r'} = g^{pqr,\; p'q'r'}_{\quad k\ell m} e^{k\ell m},
     \qquad \qquad \!
     \Delta(e^{pqr}) = f_{k\ell m,\; k'\ell' m'}^{\quad pqr} e^{k\ell m} \otimes e^{k'\ell' m'},
     \label{dualdef}
\eea
where the sum over the repeated indices is understood.
We follow \cite{FG} in order to determine the commutation relations and the coproducts
for $ x = e^{100}, y = e^{010} $ and $ z = e^{001}. $
It is not difficult to see that
\bea
  & & f_{k-1\,00,\; 100}^{\quad k\ell m} = k\delta_{\ell,0}\delta_{m,0}, \qquad\qquad\ \;
  k = 1, 2, \cdots \nn \\
  & & f_{k00\; 001}^{\ k\ell m} = \delta_{\ell,0}\delta_{m,1},  \hspace{2.6cm}
  k = 0, 1, 2, \cdots \nn \\
  & & f_{001\; 001}^{\ k\ell m} = (-1)^{k+1}2h \delta_{\ell,0}\delta_{m,0}, \qquad
  k = 1, 2, \cdots \nn
\eea
It follows that
\beq
   e^{k00} = \frac{z^k}{k!}, \qquad
   e^{k01} = \frac{x^kz}{k!}, \qquad z^2 = 2h (1-e^{-x}). \label{z2}
\eeq
We also see that
\bea
  & & f_{100,\; 010}^{\ \ 110} = f_{010,\; 100}^{\ \ 110} = 1,
      \qquad
      f_{100,\; 010}^{\ \ k00} = -f_{010,\; 100}^{\ \ k00} = (-1)^k h,
      \qquad k = 1, 2, \cdots (-1)^k, \nn \\
  & & f_{100,\; 010}^{\  k\ell m} = f_{010,\; 100}^{\  k\ell m} = 0, \qquad {\rm otherwise} \nn \\
  & & f_{100,\; 001}^{\  101} = f_{001,\; 100}^{\ 101} = 1, \qquad
      f_{100,\; 001}^{\  k\ell m} = f_{100,\; 001}^{\  k\ell m} = 0, \qquad {\rm otherwise} \nn
\eea
and
\bea
  & & f_{010,\; 001}^{\ \ 011} = f_{001,\; 010}^{\ \ 011} = 1, \qquad
      f_{010,\; 001}^{\ \ k01} = - f_{001,\; 010}^{\ \ k01} = (-1)^k\frac{h}{2},
      \qquad k = 0, 1, 2, \cdots
      \nn \\
  & & f_{010,\; 001}^{\  k\ell m} = f_{001,\; 010}^{\  k\ell m} = 0, \qquad {\rm otherwise}
      \nn
\eea
It follows the commutation relations
\beq
 [x,\; y] = 2h(e^{-x}-1), \qquad [x,\; z] = 0, \qquad [y,\; z] = hxe^{-x}.
 \label{comBs}
\eeq

  The coproducts for $ x, y $ and $ z$ are obtained from the observation
\bea
 & & g^{100\; 000}_{\ 100} = g^{000\; 100}_{\ 100} = 1, \qquad
     g^{pqr\; p'q'r'}_{\ 100} = 0, \quad {\rm otherwise} \nn \\
 & & g^{000\; 010}_{\ 010} = 1, \qquad
     g^{010\; k00}_{\ 010} = (-1)^k, \qquad k = 0, 1, 2, \cdots \nn \\
 & & g^{001\; k01}_{\ 010} = \frac{1}{4}\left(-\frac{1}{2}\right)^k, \qquad
     k = 0, 1, 2, \cdots \nn \\
 & & g^{pqr\; p'q'r'}_{\ 010} = 0, \quad {\rm otherwise} \nn \\
 & & g^{000\; 001}_{001} = 1, \qquad
     g^{001\; k00}_{001} = \left(-\frac{1}{2}\right)^k, \qquad
     k= 0, 1, 2, \cdots \nn \\
 & & g^{pqr\; p'q'r'}_{\ 001} = 0, \quad {\rm otherwise} \nn
\eea
The coproducts read
\beq
  \Delta(x) = x \otimes 1 + 1 \otimes x, \qquad
  \Delta(y) = y \otimes e^{-x} + 1 \otimes y + \frac{1}{4} z \otimes z e^{-x/2}, \qquad
  \Delta(z) = z \otimes e^{-x/2} + 1 \otimes z.
  \label{coproBs}
\eeq
We have a unique primitive element. Note that $ x, y$ are even and $ z $ are odd.
The counit and antipode are calculated from the coproducts
\bea
 & & \epsilon(x) = \epsilon(y) = \epsilon(z) = 0, \nn \\
 & & S(x) = -x, \qquad S(y) = -ye^x + \frac{h}{2}(e^x - 1), \qquad
     S(z) = -z e^{x/2}.
     \label{eSBs}
\eea
Thus $ {\cal B}^* $ is a Hopf algebra defined by the relations (\ref{z2}), (\ref{comBs}),
(\ref{coproBs}) and (\ref{eSBs}).
It is easy to see that the algebra $ {\cal B} $ is isomorphic to $ {\cal B}^*.$
The isomorphism $\rho$ is given by
\beq
  \rho(H) = \frac{1}{2h} y e^{x/2}, \qquad \rho(X) = \frac{1}{2h} x, \qquad
  \rho(V_+) = \frac{1}{4h} z e^{x/4}.
  \label{BtoBs}
\eeq

%
%%%%%%%%%%%%%%%%%%%%%%%%%%%%%%%%%%%%%%%%%%%%%%%%%%%%%%%%%%%%%%%
%
%        Concluding remarks
%
%%%%%%%%%%%%%%%%%%%%%%%%%%%%%%%%%%%%%%%%%%%%%%%%%%%%%%%%%%%%%%%%
%
\setcounter{equation}{0}
\section{Concluding remarks}

  We gave a explicit Hopf algebra structure of quantization of
$ osp(2/1) $ with the classical $r$-matrix $r_2$. It is a distinct
triangular deformation of $ osp(2/1) $ from the one by embedding
of $ sl(2) $ into $ osp(2/1) $\cite{CK}. We now have all Hopf algebra
structures obtained form the bialgebras on $ osp(2/1).$
However, one should not say that quantization of $ osp(2/1) $ is
completed, since the universal ${\cal R}$-matrix of our $ \Uh $ and
the twisting element have not been evaluated as closed form expressions.
There may be some approaches
to find the objects, for instance, power series expansion in \S 5 or
tensor operators as mentioned in \S 6. Recall also that because of the
self-duality of the quantized Borel subalgebra of $ sl(2) $ the quantum
double construction from the Borel subalgebra gives rise to the
universal $R$-matrix of whole Jordanian $ sl(2) $\cite{Vla1,Vla2,BH}.
It is observed that
the quantized Borel subalgebra ${\cal B} $ of $ osp(2/1) $ is self-dual
in \S 7. We expect that the quantum double construction gives the universal
$R$-matrix of the whole $ \Uh.$ We hope that one can write down a
closed form of the
twisting element and the universal $R$-matrix by an approach mentioned
above. This will be a future work.

%
%%%%%%%%%%%%%%%%%%%%%%%%%%%%%%%%%%%%%%%%%%%%%%%%%%%%%%%%%%%%%%%
%
%        Appendix
%
%%%%%%%%%%%%%%%%%%%%%%%%%%%%%%%%%%%%%%%%%%%%%%%%%%%%%%%%%%%%%%%%
%

\setcounter{equation}{0}
\section*{Appendix}
\renewcommand{\theequation}{A.\arabic{equation}}

In this appendix,
the fundamental and the adjoint representation matrices for $osp(2/1)$ and $ \Uh $
are summarized.

\medskip\noindent
(1) fundamental representation of $ osp(2/1) $
\beq
 \begin{array}{lll}
     J_0 = \frac{1}{2} \left(
        \begin{array}{ccc}
          1 & 0 & 0 \\
          0 & 0 & 0 \\
          0 & 0 & -1
        \end{array}
        \right)_,
      &
      v_+ = \frac{1}{2} \left(
        \begin{array}{ccc}
          0 & 1 & 0 \\
          0 & 0 & 1 \\
          0 & 0 & 0
        \end{array}
       \right)_,
      &
     v_- = \frac{1}{2} \left(
        \begin{array}{rrr}
          0 & 0 & 0 \\
          -1 & 0 & 0 \\
          0 & 1 & 0
        \end{array}
      \right)_,
     \\
     J_+ = \ \;\left(
        \begin{array}{ccc}
          0 & 0 & 1 \\
          0 & 0 & 0 \\
          0 & 0 & 0
        \end{array}
        \right)_,
     &
     J_- = \ \;\left(
        \begin{array}{ccc}
          0 & 0 & 0 \\
          0 & 0 & 0 \\
          1 & 0 & 0
        \end{array}
        \right)_.
     &
   \end{array}    \label{fundrep}
\eeq

\noindent
(2) fundamental representation of $\Uh$ by first deformation map
\beq
 \begin{array}{lll}
     H = \frac{1}{2} \left(
       \begin{array}{rrr}
          1 & 0 & h \\
          0 & 0 & 0 \\
          0 & 0 & -1
       \end{array}
       \right)_,
     &
     V_+ = \frac{1}{2} \left(
        \begin{array}{ccc}
          0 & 1 & 0 \\
          0 & 0 & 1 \\
          0 & 0 & 0
        \end{array}
        \right)_,
     &
     V_- = \frac{1}{2} \left(
        \begin{array}{rrr}
          0  & -{\textstyle \frac{h}{2}} & 0 \\
          -1 &     0                     & -{\textstyle \frac{h}{2}} \\
          0  &     1                     & 0
        \end{array}
        \right)_,
     \\
     J_+ = \ \; \left(
        \begin{array}{ccc}
          0 & 0 & 1 \\
          0 & 0 & 0 \\
          0 & 0 & 0
        \end{array}
        \right)_,
     &
     J_- = \left(
        \begin{array}{ccc}
          -{\textstyle \frac{h}{2}} & 0 & -{\textstyle \frac{h^2}{4}}\\
          0 & 0 & 0 \\
          1 & 0 & {\textstyle \frac{h}{2}}
        \end{array}
        \right)_.
      &
 \end{array}
      \label{fundrepdef}
\eeq

\noindent
(3) fundamental representation of $\Uh$ by second deformation map

  The second deformation map does not deform the fundamental
representation. The matrices for $ H, V_{\pm} $ and $ J_{\pm} $
are the same as (\ref{fundrep}).

\medskip\noindent
(4) adjoint representation of $ ops(2/1) $

  The representation basis is the algebra itself. The basis
is taken as $ (J_+, v_+, J_0, v_-, J_-). $
\beq
 \begin{array}{ll}
  J_+ = \left(
    \begin{array}{rrrrr}
      0 & 0 & -1 & 0 & 0 \\
      0 & 0 &  0 & 1 & 0 \\
      0 & 0 &  0 & 0 & 2 \\
      0 & 0 &  0 & 0 & 0 \\
      0 & 0 &  0 & 0 & 0
    \end{array}
    \right)_,
 &
 v_+ = \left(
    \begin{array}{rrrrr}
      0 & \frac{1}{2} & 0 & 0 & 0 \\
      0 & 0 &  -\frac{1}{2} & 0 & 0 \\
      0 & 0 &  0 & -\frac{1}{2} & 0 \\
      0 & 0 &  0 & 0 & -1 \\
      0 & 0 &  0 & 0 & 0
    \end{array}
    \right)_,
  \\
 J_0 = diag\left(1, \frac{1}{2}, 0, -\frac{1}{2}, -1 \right), \label{matundef} \\
 &
 \\
 v_- = \left(
    \begin{array}{rrrrr}
      0  & 0 &  0 & 0 & 0 \\
      -1 & 0 &  0 & 0 & 0 \\
      0  & -\frac{1}{2} &  0 & 0 & 0 \\
      0  & 0 & \frac{1}{2} & 0 & 0 \\
      0  & 0 &  0 & -\frac{1}{2} & 0
    \end{array}
    \right)_,\quad
 &
 J_- = \left(
    \begin{array}{rrrrr}
      0  & 0 &  0 & 0 & 0 \\
      0  & 0 &  0 & 0 & 0 \\
      -2 & 0 &  0 & 0 & 0 \\
      0  & 1 &  0 & 0 & 0 \\
      0  & 0 &  1 & 0 & 0
    \end{array}
    \right)_.
  \end{array}
  \label{adjrep}
\eeq

\noindent
(5) adjoint representation of $\Uh$ by first deformation map
\beq
 \begin{array}{ll}
 X = \left(
    \begin{array}{rrrrr}
      0 & 0 & -1 & 0 & -2h \\
      0 & 0 &  0 & 1 & 0 \\
      0 & 0 &  0 & 0 & 2 \\
      0 & 0 &  0 & 0 & 0 \\
      0 & 0 &  0 & 0 & 0
    \end{array}
    \right)_,
 &
 V_+ = \left(
    \begin{array}{rrrrr}
      0 & \frac{1}{2} & 0 & \frac{h}{4} & 0 \\
      0 & 0 &  -\frac{1}{2} & 0 & -\frac{h}{2} \\
      0 & 0 &  0 & -\frac{1}{2} & 0 \\
      0 & 0 &  0 & 0 & -1 \\
      0 & 0 &  0 & 0 & 0
    \end{array}
    \right)_,
 \\
 H = \left(
    \begin{array}{rrrrr}
     1 & 0 & 0 & 0 & -h^2 \\
     0 & \frac{1}{2} & 0 & \frac{h}{2} & 0 \\
     0 & 0 & 0 & 0 & 2h \\
     0 & 0 & 0 & - \frac{1}{2} & 0 \\
     0 & 0 & 0 & 0 & -1
    \end{array}
    \right)_,
 &
 V_- = \left(
    \begin{array}{rrrrr}
      0  & 0 &  0 & -\frac{5h^2}{16} & 0 \\
      -1 & 0 &  -\frac{h}{4} & 0 & -\frac{h^2}{8} \\
      0  & -\frac{1}{2} &  0 & \frac{3h}{4} & 0 \\
      0  & 0 & \frac{1}{2} & 0 & h \\
      0  & 0 &  0 & -\frac{1}{2} & 0
    \end{array}
    \right)_,  \nn \\
 \\
 \multicolumn{2}{l}{
 Y = \left(
    \begin{array}{rrrrr}
      0  & 0 &  \frac{5h^2}{8} & 0 & \frac{5h^3}{4} \\
      0  & -\frac{h}{2} &  0 & -\frac{3h^2}{4} & 0 \\
      -2 & 0 &  -2h & 0 & -\frac{13h^2}{4} \\
      0  & 1 &  0 & \frac{h}{2} & 0 \\
      0  & 0 &  1 & 0 & 2h
    \end{array}
    \right)_.
 }
 \end{array} \label{1stdefmat}
\eeq

\noindent
(6) adjoint representation of $\Uh$ by second deformation map
\beq
 \begin{array}{ll}
 X = \left(
    \begin{array}{rrrrr}
      0 & 0 & -1 & 0 & 0 \\
      0 & 0 &  0 & 1 & 0 \\
      0 & 0 &  0 & 0 & 2 \\
      0 & 0 &  0 & 0 & 0 \\
      0 & 0 &  0 & 0 & 0
    \end{array}
    \right)_,
 &
 V_+ = \left(
    \begin{array}{rrrrr}
      0 & \frac{1}{2} & 0 & 0 & 0 \\
      0 & 0 &  -\frac{1}{2} & 0 & 0 \\
      0 & 0 &  0 & -\frac{1}{2} & 0 \\
      0 & 0 &  0 & 0 & -1 \\
      0 & 0 &  0 & 0 & 0
    \end{array}
    \right)_,
 \\
 H = diag\left(1, \frac{1}{2}, 0, -\frac{1}{2}, -1 \right),
 &
 \\
 V_- = \left(
    \begin{array}{rrrrr}
      0  & 0 &  0 & -\frac{h^2}{16} & 0 \\
      -1 & 0 &  0 & 0 & -\frac{h^2}{8} \\
      0  & -\frac{1}{2} &  0 & 0 & 0 \\
      0  & 0 & \frac{1}{2} & 0 & 0 \\
      0  & 0 &  0 & -\frac{1}{2} & 0
    \end{array}
    \right)_,
 &
 Y = \left(
    \begin{array}{rrrrr}
      0  & 0 &  \frac{h^2}{8} & 0 & 0 \\
      0  & 0 &  0 & -\frac{h^2}{2} & 0 \\
      -2 & 0 &  0 & 0 & -\frac{h^2}{4} \\
      0  & 1 &  0 & 0 & 0 \\
      0  & 0 &  1 & 0 & 0
    \end{array}
    \right)_.
 \end{array}
 \label{2nddefmat}
\eeq
$X, V_+$ and $ H $ remain undeformed by this map.
%
%%%%%%%%%%%%%%%%%%%%%%%%%%%%%%%%%%%%%%%%%%%%%%%%%%%%%%%%%%%%%%%%
%
%  References
%
%%%%%%%%%%%%%%%%%%%%%%%%%%%%%%%%%%%%%%%%%%%%%%%%%%%%%%%%%%%%%%%%
%

\end{document}